\documentclass[a4paper]{amsart}  

\usepackage{amssymb} \usepackage{amscd} \usepackage{times}

\def\zbb{\mathbb{Z}}  
  
  \def\phi{\varphi}
 \def\p1{{\mathbb{P}^1_\zbb}}

\newtheorem{Theorem}{\quad Theorem}[section]

\newtheorem{Corollary}[Theorem]{\quad Corollary}
\newtheorem{Lemma}[Theorem]{\quad Lemma}
\newtheorem{Proposition}[Theorem]{\quad Proposition}

\newcommand{\be} {\begin{equation}}
\newcommand{\ee} {\end{equation}}

\begin{document}

\title{ Harnack type inequality for an elliptic equation.}

\author{Samy Skander Bahoura}

\address{Department of Mathematics, Pierre et Marie Curie University, 75005 Paris France.}
              
\email{samybahoura@yahoo.fr, bahoura@ccr.jussieu.fr} 

\date{}

\maketitle

\begin{abstract}

We give a $ \sup \times \inf $ inequality for an elliptic equation.

\end{abstract}

\section{Introduction and Main Results}

We are on Riemannian manifold $ (M,g) $ of dimension $ n \geq 3 $. In this paper we denote $ \Delta_g = -\nabla^j(\nabla_j) $ the Laplace-Beltrami operator and $ N=\frac{2n}{n-2} $.

We consider the following equation

\be \Delta_g u=V u^{N-1}+u^{\alpha},\,\, u >0 . \ee

Where $ V $ is a function and $ \alpha \in ]\frac{n}{n-2}, \frac{n+2}{n-2}[ $.

For $ a, b, A >0 $, we consider a sequence $ (u_i, V_i)_i $ of solutions of the previous equation with the following conditions:

$$ 0 <a \leq V_i \leq b < +\infty, $$

$$ ||\nabla V_i||_{\infty} \leq A. $$

Here we study some properties of this nonlinear elliptic equation. We try to find some estimates of type $ \sup \times \inf $. We denote by $ S_g $ the scalar curvature.

There are many existence and compactness results which concern this type of equations, see for example [1-21]. In particluar in [1], we can find some results about the Yamabe equation and the Prescribed scalar curvature equation. Many methods where used to solve these problems, as a variationnal approach and some other topological methods. Note that the problems come from the nonlinearity of the critical Sobolev exponent. We can find in [1] some uniform estimates for various equations on the unit sphere or for the Monge-Ampere equation. Note that Tian and Siu proved uniform upper and lower bounds for the $ \sup+\inf $ for the Monge-Ampere equation under some condition on the Chern class, see [1]. In the case of the Scalar curvature equation and in dimension 2 Shafrir used the isoperimetric inequality of Alexandrov to prove an inequality of type $ \sup+\inf $ with only $ L^{\infty} $ assumption on the prescribed curvature, see [21]. The result of Shafrir is an extention of a result of Brezis and Merle, see [4] and later, Brezis-Li-Shafrir proved a sharp $ \sup +\inf  $ inequality  for the same equation with Lipschitzian assumption on the prescribed scalar curvature, see [3].  Li in[17] extend the previous last result to compact Riemannian surfaces. In the higher dimensional case, we can find in [15] a proof of the $  \sup \times \inf  $ inequality in the constant case for the scalar curvature equation on open set of $ {\mathbb R}^n $. We have various estimates in [2] when we consider the nonconstant case. To prove our result, we use a blow-up analysis and the moving-plane method, based on the maximum principle and the Hopf Lemma as showed in [2, 3, 15, 17], and a condition on the scalar curvature is sufficient to prove the estimate.

Our main result is:

\begin{Theorem} Assume $  S_g >0 $ on $ M $, then, for every compact $ K $ of $ M $, there exist a positive constant $ c=c(\alpha, a, b, A, K,M,n,g) $ such that:

$$ \sup_K u_i \times \inf_M u_i \leq c. $$
\end{Theorem}

{\bf Remark:} in the case where $ (M,g)=(\Omega \subset {\mathbb R}^n, \delta) $ an open set of the euclidean space with the flat metric,  we have the same inequality on compact sets of $ \Omega $ in this case the scalar curvature $ S_{\delta} \equiv 0 $, see [2]. 

\bigskip

If we consider the Green function $ G $ of the Laplacian with Dirichlet condition on small balls of $ M $, we can have a positive lower bound for $ G $ and we have the following corollary:

\begin{Corollary} Assume $  S_g >0 $ on $ M $, then, for every compact $ K $ of $ M $, there exist a positive constant $ c'=c'(\alpha, a, b, A, K,M,n,g) $ such that:

$$ \int_K u_i^{\frac{2n}{n-2}} dv_g \leq c'. $$

\end{Corollary}

\section{Proof of the theorem.}

Let us consider $ x_0 \in M $, by a conformal change of the metric $ \tilde g = \phi^{4/(n-2)} g $ with $ \phi >0 $ we can consider an equation of type:

\be \Delta_{\tilde g}u + R_{\tilde g} u=V u^{N-1}+ \phi^{\alpha+1-N}u^{\alpha}+ R_g \phi^{2-N} u,\,\, u >0 . \ee 

with,

$$ Ricci_{\tilde g}(x_0)= 0. $$

Here; $ R_g= \dfrac{n-2}{4(n-1)} S_g $ and  $ R_{\tilde g}= \dfrac{n-2}{4(n-1)} S_{\tilde g} $

It is clear see the computations in a previous paper [2], it is sufficient to consider an equation of type:

\be \Delta_{\tilde g}u =V u^{N-1}+ u^{\alpha}+ \mu u,\,\, u >0 . \ee

with,

$ Ricci \equiv 0 $ and $ \mu > 0 $.

\bigskip

{\underbar {\bf Part I: The metric in polar coordinates.}}

\bigskip

Let $ (M,g) $ a Riemannian manifold. We note $ g_{x,ij} $ the local expression of the metric $ g $  in the exponential map centered in $ x $.

We are concerning by the polar coordinates expression of the metric. Using Gauss lemma, we can write:

$$ g=ds^2=dt^2+g_{ij}^k(r,\theta)d\theta^id\theta^j=dt^2+r^2{\tilde g}_{ij}^k(r,\theta)d\theta^id\theta^j=g_{x,ij}dx^idx^j, $$

in a polar chart with origin $ x $", $ ]0,\epsilon_0[\times U^k $, with $ ( U^k, \psi) $ a chart of $ {\mathbb S}_{n-1} $. We can write the element volume:

$$ dV_g=r^{n-1}\sqrt {|{\tilde g}^k|}dr d \theta^1 \ldots d \theta^{n-1} = \sqrt {[det(g_{x,ij})]}dx^1 \ldots dx^n, $$

then,

$$ dV_g=r^{n-1} \sqrt { [det(g_{x,ij})]}[\exp_x(r\theta)]\alpha^k(\theta)dr d\theta^1 \ldots d \theta^{n-1} , $$

where, $ \alpha^k $ is such that, $ d\sigma_{{\mathbb S}_{n-1}}=\alpha^k(\theta) d \theta^1 \ldots d \theta^{n-1} . $ (Riemannian volume element of the sphere in the chart $ (U^k,\psi) $ ).

\bigskip

Then,

$$ \sqrt { |{\tilde g}^k|}=\alpha^k(\theta) \sqrt {[det(g_{x,ij})]}. $$

Clearly, we have the following proposition:

\bigskip

\begin{Proposition} Let $ x_0 \in M $, there exist $ \epsilon_1>0 $ and if we reduce $ U^k $, we have:

$$ |\partial_r{\tilde g}_{ij}^k(x,r,\theta)|+|\partial_r\partial_{\theta^m}{\tilde g}_{ij}^k(x,r, \theta)| \leq C r,\,\, \forall \,\, x\in B(x_0,\epsilon_1) \,\, \forall \,\, r\in [0,\epsilon_1], \,\, \forall \,\, \theta \in U^k.$$

and,

$$ |\partial_r|{\tilde g}^k|(x,r,\theta)|+\partial_r \partial_{\theta^m} |{\tilde g}^k|(x,r,\theta)\leq C r,\,\, \forall \,\, x\in B(x_0,\epsilon_1) \,\, \forall \,\, r\in [0,\epsilon_1], \,\, \forall \,\, \theta \in U^k. $$
\end{Proposition}

\underbar {\bf Remark:} 

\bigskip

$ \partial_r [ \log \sqrt { |{\tilde g}^k |}] $ is a local function of $ \theta $, and the restriction of the global function on the sphere $ {\mathbb S}_{n-1} $, $ \partial_r [\log \sqrt { det(g_{x, ij})}] $. We will note, $ J(x,r,\theta)=\sqrt { det(g_{x, ij})} $.

\bigskip

{\underbar { \bf Part II: The laplacian in polar coordinates}}

\bigskip

Let's write the laplacian in $ [0,\epsilon_1]\times U^k $,

$$ -\Delta = \partial_{rr}+\dfrac{n-1}{r}\partial_r+ \partial_r [\log \sqrt { |{\tilde g^k|}] }\partial_r+\dfrac{1}{r^2 \sqrt {|{\tilde g}^k|}}\partial_{\theta^i}({\tilde g}^{\theta^i \theta^j}\sqrt { |{\tilde g}^k|}\partial_{\theta^j}) . $$

We have,

$$ -\Delta = \partial_{rr}+\dfrac{n-1}{r}\partial_r+ \partial_r \log J(x,r,\theta)\partial_r+ \dfrac{1}{r^2 \sqrt {|{\tilde g}^k|}}\partial_{\theta^i}({\tilde g}^{\theta^i \theta^j}\sqrt { |{\tilde g}^k|}\partial_{\theta^j}) . $$

We write the laplacian ( radial and angular decomposition),

$$ -\Delta = \partial_{rr}+\dfrac{n-1}{r} \partial_r+\partial_r [\log J(x,r,\theta)] \partial_r-\Delta_{{S}_r(x)}, $$

where $ \Delta_{ S_r(x)} $ is the laplacian on the sphere $ {S}_r(x) $. 

\bigskip

We set $ L_{\theta}(x,r)(...)=r^2\Delta_{ S_r(x)}(...)[\exp_x(r\theta)] $, clearly, this operator is a laplacian on $ {\mathbb S}_{n-1} $ for particular metric. We write,

$$ L_{\theta}(x,r)=\Delta_{g_{x,r, {}_{{\mathbb S}_{n-1}}}}, $$

and,

$$ \Delta = \partial_{rr}+\dfrac{n-1}{r} \partial_r+\partial_r [ J(x,r,\theta)] \partial_r - \dfrac{1}{r^2} L_{\theta}(x,r) . $$

If, $ u $ is function on $ M $, then, $ \bar u(r,\theta)=u[\exp_x(r\theta)] $ is the corresponding function in polar coordinates centered in $ x $. We have,

\be -\Delta u =\partial_{rr} \bar u+\dfrac{n-1}{r} \partial_r \bar u+\partial_r [ J(x,r,\theta)] \partial_r \bar u-\dfrac{1}{r^2}L_{\theta}(x,r)\bar u .\ee

{\underbar { \bf Part III: "Blow-up" and "Moving-plane" methods }} 

\bigskip

\underbar {\bf The "blow-up" analysis}

\bigskip

Let, $ (u_i)_i $ a sequence of functions on $ M $ such that,

\be \Delta_{\tilde g}u_i + R_{\tilde g} u_i=V_i u_i^{N-1}+ \phi^{\alpha-N}u_i^{\alpha}+ R_g \phi^{2-N} u_i,\,\, u_i>0,\,\, N=\dfrac{2n}{n-2}, \ee

It is sufficient to consider an equation of type:

\be \Delta_{\tilde g}u =V u^{N-1}+ u^{\alpha}+ \mu u,\,\, u >0 . \ee

with $ Ricci \equiv 0 $ and $ \mu >0 $.

We argue by contradiction and we suppose that $ \sup \times \inf $ is not bounded.

\smallskip

We assume that:

\bigskip

$ \forall \,\, c,R >0 \,\, \exists \,\, u_{c,R} $ solution of $ (E) $ such that:

$$ R^{n-2} \sup_{B(x_0,R)} u_{c,R} \times \inf_M u_{c,R} \geq c. \qquad (H) $$

\bigskip

\begin{Proposition}

There exist a sequence of points $ (y_i)_i $, $ y_i \to x_0 $ and two sequences of positive real number $ (l_i)_i, (L_i)_i $, $ l_i \to 0 $, $ L_i \to +\infty $, such that if we consider $ v_i(y)=\dfrac{u_i[\exp_{y_i}(y)]}{u_i(y_i)} $, we have:

$$ i) \qquad 0 < v_i(y) \leq  \beta_i \leq 2^{(n-2)/2}, \,\, \beta_i \to 1. $$

$$ ii) \qquad v_i(y)  \to \left ( \dfrac{1}{1+{|y|^2}} \right )^{(n-2)/2}, \,\, {\rm uniformly \,\, on\,\, every \,\, compact \,\, set \,\, of } \,\, {\mathbb R}^n . $$

$$ iii) \qquad l_i^{(n-2)/2} [u_i(y_i)] \times \inf_M u_i \to +\infty $$
\end{Proposition}

\underbar {\bf Proof:}

Without loss of generality, we can assume that:

$$ V(x_0)=n(n-2). $$

We use the hypothesis $ (H) $. We can take two sequences $ R_i>0, R_i \to 0 $ and $ c_i \to +\infty $, such that,

$$ {R_i}^{(n-2)} \sup_{B(x_0,R_i)} u_i \times \inf_M u_i \geq c_i \to +\infty. $$

Let, $ x_i \in  { B(x_0,R_i)} $, such that $ \sup_{B(x_0,R_i)} u_i=u_i(x_i) $ and $ s_i(x)=[R_i-d(x,x_i)]^{(n-2)/2} u_i(x), x\in B(x_i, R_i) $. Then, $ x_i \to x_0 $.

\bigskip

We have, 

$$ \max_{B(x_i,R_i)} s_i(x)=s_i(y_i) \geq s_i(x_i)={R_i}^{(n-2)/2} u_i(x_i)\geq \sqrt {c_i}  \to + \infty. $$ 

\bigskip

Set :

$$ l_i=R_i-d(y_i,x_i),\,\, \bar u_i(y)= u_i [\exp_{y_i}(y)],\,\,  v_i(z)=\dfrac{u_i [ \exp_{y_i}\left ( z/[u_i(y_i)]^{2/(n-2)} \right )] } {u_i(y_i)}. $$

Clearly, $ y_i \to x_0 $. We obtain:

$$ L_i= \dfrac{l_i}{(c_i)^{1/2(n-2)}} [u_i(y_i)]^{2/(n-2)}=\dfrac{[s_i(y_i)]^{2/(n-2)}}{c_i^{1/2(n-2)}}\geq \dfrac{c_i^{1/(n-2)}}{c_i^{1/2(n-2)}}=c_i^{1/2(n-2)}\to +\infty. $$

\bigskip

If $ |z|\leq L_i $, then $ y=\exp_{y_i}[z/ [u_i(y_i)]^{2/(n-2)}] \in B(y_i,\delta_i l_i) $ with $ \delta_i=\dfrac{1}{(c_i)^{1/2(n-2)}} $ and $ d(y,y_i) < R_i-d(y_i,x_i) $, thus, $ d(y, x_i) < R_i $ and, $ s_i(y)\leq s_i(y_i) $, we can write,

$$ u_i(y) [R_i-d(y,y_i)]^{(n-2)/2} \leq u_i(y_i) (l_i)^{(n-2)/2}. $$

But, $ d(y,y_i) \leq \delta_i l_i $, $ R_i >l_i$ and $ R_i-d(y, y_i) \geq R_i-\delta_i l_i>l_i-\delta_i l_i=l_i(1-\delta_i) $, we obtain,

$$ 0 < v_i(z)=\dfrac{u_i(y)}{u_i(y_i)} \leq \left [ \dfrac{l_i}{l_i(1-\delta_i)} \right ]^{(n-2)/2}\leq 2^{(n-2)/2} . $$

We set, $ \beta_i=\left ( \dfrac{1}{1-\delta_i} \right )^{(n-2)/2} $, clearly $ \beta_i \to 1 $.

\bigskip

The function $ v_i $ is solution of:

$$ -g^{jk}[\exp_{y_i}(y)]\partial_{jk} v_i-\partial_k \left [ g^{jk}\sqrt { |g| } \right ][\exp_{y_i}(y)]\partial_j v_i=\dfrac{1}{[u_i(y_i)]^{N-1-\alpha}} v_i^{\alpha}+\dfrac{ \mu}{[u_i(y_i)]^{4/(n-2)}} v_i+V_i{v_i}^{N-1}, $$

By elliptic estimates and Ascoli, Ladyzenskaya theorems, $ ( v_i)_i $ converge uniformely on each compact to the function $ v $ solution on $ {\mathbb R}^n $ of, 

\be \Delta v=n(n-2)v^{N-1}, \,\, v(0)=1,\,\, 0 \leq v\leq 1\leq 2^{(n-2)/2}, \ee

By using maximum principle, we have $ v>0 $ on $ {\mathbb R}^n $, the result of Caffarelli-Gidas-Spruck ( see [6]) give, $ v(y)=\left ( \dfrac{1}{1+{|y|^2}} \right )^{(n-2)/2} $. We have the same properties for $ v_i $ in the previous paper [2].

\bigskip

\underbar {\bf Polar coordinates and "moving-plane" method}

\bigskip

Let, 

$$ w_i(t,\theta)=e^{(n-2)/2}\bar u_i(e^t,\theta) = e^{(n-2)t/2}u_io\exp_{y_i}(e^t\theta), \,\, {\rm et} \,\, a(y_i,t,\theta)=\log J(y_i,e^t,\theta). $$ 

We set $ \delta = \dfrac{(n+2)-(n-2)\alpha}{2} $.

\smallskip

\begin{Lemma}

The function $ w_i $ is solution of:

\be  -\partial_{tt} w_i-\partial_t a \partial_t w_i-L_{\theta}(y_i,e^t)+c w_i=V_iw_i^{N-1}+e^{\delta t} w_i^{\alpha}+\mu e^{2t}w_i, \ee
 
with,

 $$ c = c(y_i,t,\theta)=\left ( \dfrac{n-2}{2} \right )^2+ \dfrac{n-2}{2} \partial_t a. $$ 

\end{Lemma}

\underbar {\bf Proof:}

\smallskip

We write:

$$ \partial_t w_i=e^{nt/2}\partial_r \bar u_i+\dfrac{n-2}{2} w_i,\,\, \partial_{tt} w_i=e^{(n+2)t/2} \left [\partial_{rr} \bar u_i+\dfrac{n-1}{e^t}\partial_r \bar u_i \right ]+\left ( \dfrac{n-2}{2} \right )^2 w_i. $$

$$ \partial_t a =e^t\partial_r \log J(y_i,e^t,\theta), \partial_t a \partial_t w_i=e^{(n+2)t/2}\left [ \partial_r \log J\partial_r \bar u_i \right ]+\dfrac{n-2}{2} \partial_t a w_i.$$

the lemma  is proved.

\bigskip

Now we have, $ \partial_t a=\dfrac{ \partial_t b_1}{b_1} $, $ b_1(y_i,t,\theta)=J(y_i,e^t,\theta)>0 $,

\bigskip

We can write,

$$ -\dfrac{1}{\sqrt {b_1}}\partial_{tt} (\sqrt { b_1} w_i)-L_{\theta}(y_i,e^t)w_i+[c(t)+ b_1^{-1/2} b_2(t,\theta)]w_i=V_iw_i^{N-1}+e^{\delta t} w_i^{\alpha}+\mu e^{2t}w_i, $$

where, $ b_2(t,\theta)=\partial_{tt} (\sqrt {b_1})=\dfrac{1}{2 \sqrt { b_1}}\partial_{tt}b_1-\dfrac{1}{4(b_1)^{3/2}}(\partial_t b_1)^2 .$

\bigskip

Let,

$$ \tilde w_i=\sqrt {b_1} w_i, $$

\begin{Lemma}

The function $ \tilde w_i $ is solution of:

$$ -\partial_{tt} \tilde w_i+\Delta_{g_{y_i, e^t, {}_{{\mathbb S}_{n-1}}}} (\tilde w_i)+2\nabla_{\theta}(\tilde  w_i) .\nabla_{\theta} \log (\sqrt {b_1})+(c+b_1^{-1/2} b_2-c_2) \tilde w_i= $$

\be = V_i\left (\dfrac{1}{b_1} \right )^{(N-2)/2} {\tilde w_i}^{N-1}+e^{\delta t}\left (\dfrac{1}{b_1} \right )^{(\alpha-1)/2} {\tilde w_i}^{\alpha}+\mu e^{2t}{\tilde w_i}, \ee

where, $ c_2 =[\dfrac{1}{\sqrt {b_1}} \Delta_{g_{y_i, e^t, {}_{{\mathbb S}_{n-1}}}}(\sqrt{b_1}) + |\nabla_{\theta} \log (\sqrt {b_1})|^2] . $

\end{Lemma}

\underbar {\bf Proof:}

\bigskip

We have: 

$$ -\partial_{tt} \tilde w_i-\sqrt {b_1} \Delta_{g_{y_i, e^t, {}_{{\mathbb S}_{n-1}}}} w_i+(c+b_2) \tilde w_i= V_i\left (\dfrac{1}{b_1} \right )^{(N-2)/2} {\tilde w_i}^{N-1}+$$

$$ +e^{\delta t}\left (\dfrac{1}{b_1} \right )^{(\alpha-1)/2} {\tilde w_i}^{\alpha}+\mu e^{2t}{\tilde w_i}, $$

But,

$$ \Delta_{g_{y_i, e^t, {}_{{\mathbb S}_{n-1}}}} (\sqrt {b_1} w_i)=\sqrt {b_1} \Delta_{g_{y_i, e^t, {}_{{\mathbb S}_{n-1}}}} w_i-2 \nabla_{\theta} w_i .\nabla_{\theta} \sqrt {b_1}+ w_i \Delta_{g_{y_i, e^t, {}_{{\mathbb S}_{n-1}}}}(\sqrt {b_1}), $$

and,

$$ \nabla_{\theta} (\sqrt {b_1} w_i)=w_i \nabla_{\theta} \sqrt {b_1}+ \sqrt {b_1} \nabla_{\theta} w_i, $$

we deduce than,

$$  \sqrt {b_1} \Delta_{g_{y_i, e^t, {}_{{\mathbb S}_{n-1}}}} w_i= \Delta_{g_{y_i, e^t, {}_{{\mathbb S}_{n-1}}}} (\tilde w_i)+2\nabla_{\theta}(\tilde  w_i) .\nabla_{\theta} \log (\sqrt {b_1})-c_2 \tilde w_i, $$

with $ c_2=[\dfrac{1}{\sqrt {b_1}} \Delta_{g_{y_i, e^t, {}_{{\mathbb S}_{n-1}}}}(\sqrt{b_1}) + |\nabla_{\theta} \log (\sqrt {b_1})|^2] . $ The lemma is proved.

\bigskip

\underbar {\bf The "moving-plane" method:}

\bigskip

Let $ \xi_i $ a real number,  and suppose $ \xi_i \leq t $. We set $ t^{\xi_i}=2\xi_i-t $ and $ \tilde w_i^{\xi_i}(t,\theta)=\tilde w_i(t^{\xi_i},\theta) $.

\bigskip

We have, 

$$ -\partial_{tt} \tilde w_i^{\xi_i}+\Delta_{g_{y_i, e^{t^{\xi_i}} {}_{{\mathbb S}_{n-1}}}} (\tilde w_i)+2\nabla_{\theta}(\tilde  w_i^{\xi_i}) .\nabla_{\theta} \log (\sqrt {b_1}) \tilde w_i^{\xi_i}+[c(t^{\xi_i})+b_1^{-1/2}(t^{\xi_i},.)b_2(t^{\xi_i})-c_2^{\xi_i}] \tilde w_i^{\xi_i}= $$

$$ =V_i\left (\dfrac{1}{b_1^{\xi_i}} \right )^{(N-2)/2} { ({\tilde w_i}^{\xi_i}) }^{N-1}+ $$

$$ +e^{\delta t}\left (\dfrac{1}{b_1} \right )^{(\alpha-1)/2} {\tilde w_i}^{\alpha}+\mu e^{2t}{\tilde w_i}, $$

By using the same arguments than in [2], we have:

\bigskip

\begin{Proposition}

We have:

$$ 1)\,\,\, \tilde w_i(\lambda_i,\theta)-\tilde w_i(\lambda_i+4,\theta) \geq \tilde k>0, \,\, \forall \,\, \theta \in {\mathbb S}_{n-1}. $$

For all $ \beta >0 $, there exist $ c_{\beta} >0 $ such that:

$$ 2) \,\,\, \dfrac{1}{c_{\beta}} e^{(n-2)t/2}\leq \tilde w_i(\lambda_i+t,\theta) \leq c_{\beta}e^{(n-2)t/2}, \,\, \forall \,\, t\leq \beta, \,\, \forall \,\, \theta \in {\mathbb S}_{n-1}. $$

\end{Proposition}

We set,

$$ \bar Z_i=-\partial_{tt} (...)+\Delta_{g_{y_i, e^t, {}_{{\mathbb S}_{n-1}}}} (...)+2\nabla_{\theta}(...) .\nabla_{\theta} \log (\sqrt {b_1})+(c+b_1^{-1/2} b_2-c_2)(...) $$

{\bf Remark:} In the operator $ \bar Z_i $, by using the proposition 3, the coeficient $ c+b_1^{-1/2}b_2-c_2 $ satisfies:

$$ c+b_1^{-1/2}b_2-c_2 \geq k'>0,\,\, {\rm pour }\,\, t<<0, $$

it is fundamental if we want to apply Hopf maximum principle.

We set $ \delta = \dfrac{(n+2)-(n-2)\alpha}{2} $.

\bigskip

\underbar {\bf Goal:}

\bigskip

Like in [2], we have elliptic second order operator. Here it is $ \bar Z_i $, the goal is to use the "moving-plane" method to have a contradiction. For this, we must have:

$$ \bar Z_i(\tilde w_i^{\xi_i}-\tilde w_i) \leq 0, \,\, {\rm if} \,\, \tilde w_i^{\xi_i}-\tilde w_i \leq 0. $$

We write:

$$ \bar Z_i(\tilde w_i^{\xi_i}-\tilde w_i)= (\Delta_{g_{y_i, e^{t^{\xi_i}}, {}_{{\mathbb S}_{n-1}}}}-\Delta_{g_{y_i, e^{t}, {}_{{\mathbb S}_{n-1}}}}) (\tilde w_i^{\xi_i})+ $$

$$ +2(\nabla_{\theta, e^{t^{\xi_i}}}-\nabla_{\theta, e^t})(w_i^{\xi_i}) .\nabla_{\theta, e^{t^{\xi_i}}} \log (\sqrt {b_1^{\xi_i}})+ 2\nabla_{\theta,e^t}(\tilde w_i^{\xi_i}).\nabla_{\theta, e^{t^{\xi_i}}}[\log (\sqrt {b_1^{\xi_i}})-\log \sqrt {b_1}]+ $$ 

$$ +2\nabla_{\theta,e^t} w_i^{\xi_i}.(\nabla_{\theta,e^{t^{\xi_i}}}-\nabla_{\theta,e^t})\log \sqrt {b_1}- [(c+b_1^{-1/2} b_2-c_2)^{\xi_i}-(c+b_1^{-1/2}b_2-c_2)]\tilde w_i^{\xi_i} + $$

$$ + V_i^{\xi_i}\left ( \dfrac{1}{b_1^{\xi_i}} \right )^{(N-2)/2} ({\tilde w_i}^{\xi_i})^{N-1}-V_i\left ( \dfrac{1}{b_1} \right )^{(N-2)/2} {\tilde w_i}^{N-1}+ $$
$$ +e^{\delta {t^{\xi_i}}} {b_1^{\xi_i}}^{(1-\alpha)/2}({\tilde w_i}^{\xi_i})^{\alpha}-e^{\delta t} b_1^{(1-\alpha)/2}({\tilde w_i})^{\alpha}+\mu( e^{2t^{\xi_i}}{\tilde w_i}^{\xi_i}-e^{2t}{\tilde w_i}) \qquad (***1) $$

Clearly, we have:

\begin{Lemma}

$$ b_1(y_i,t,\theta)=1-\dfrac{1}{3} Ricci_{y_i}(\theta,\theta)e^{2t}+\ldots, $$

$$ R_g(e^t\theta)=R_g(y_i) + <\nabla R_g(y_i)|\theta > e^t+\dots . $$
\end{Lemma}

According to proposition 1 and lemma 3,

\bigskip

\begin{Proposition}

$$ \bar Z_i(\tilde w_i^{\xi_i}-\tilde w_i) \leq {\tilde A}(e^{t}-e^{t^{\xi}})(\tilde w_i^{\xi_i})^{N-1})+(1/2)(e^{\delta t^{\xi_i}}-e^{\delta t})(\tilde w_i^{\xi_i})^{\alpha}+ $$

$$ +C|e^{2t}-e^{2t^{\xi_i}}|\left [|\nabla_{\theta} {\tilde w_i}^{\xi_i}| + |\nabla_{\theta}^2(\tilde w_i^{\xi_i})|+ o(1)[\tilde w_i^{\xi_i}+(\tilde w_i^{\xi_i})^{N-1}+ (\tilde w_i^{\xi_i})^{\alpha}] +\mu \tilde w_i^{\xi_i} \right ]. $$
\end{Proposition}

\underbar {\bf Proof:}

\bigskip

We use proposition 1, we have:

$$ a(y_i,t,\theta)=\log J(y_i,e^t,\theta)=\log b_1, |\partial_t b_1(t)|+|\partial_{tt} b_1(t)|+|\partial_{tt} a(t)|\leq C e^{2t}, $$

and,

$$ |\partial_{\theta_j} b_1|+|\partial_{\theta_j,\theta_k} b_1|+\partial_{t,\theta_j}b_1|+|\partial_{t,\theta_j,\theta_k} b_1|\leq C e^{2t}, $$

then,

$$ |\partial_t b_1(t^{\xi_i})-\partial_t b_1(t)|\leq C'|e^{2t}-e^{2t^{\xi_i}}|,\,\, {\rm on} \,\, ]-\infty, \log \epsilon_1]\times {\mathbb S}_{n-1},\forall \,\, x\in B(x_0,\epsilon_1) $$

Locally,

$$ \Delta_{g_{y_i, e^t, {}_{{\mathbb S}_{n-1}}}}= L_{\theta}(y_i,e^t)=-\dfrac{1}{\sqrt {|{\tilde g}^k(e^t,\theta)|}}\partial_{\theta^l}[{\tilde g}^{\theta^l \theta^j}(e^t,\theta)\sqrt { |{\tilde g}^k(e^t,\theta)|}\partial_{\theta^j}] . $$

Thus, in $ [0,\epsilon_1]\times U^k $, we have,

$$ A_i=\left [{ \left [ \dfrac{1}{\sqrt {|{\tilde g}^k|}}\partial_{\theta^l}({\tilde g}^{\theta^l \theta^j}\sqrt { |{\tilde g}^k|}\partial_{\theta^j}) \right ] }^{\xi_i}- \dfrac{1}{\sqrt {|{\tilde g}^k|}}\partial_{\theta^l}({\tilde g}^{\theta^l \theta^j}\sqrt { |{\tilde g}^k|}\partial_{\theta^j}) \right ](\tilde w_i^{\xi_i}) $$

then, $ A_i=B_i+D_i $ with,

$$ B_i=\left [ {\tilde g}^{\theta^l \theta^j}(e^{t^{\xi_i}}, \theta)-{\tilde g}^{\theta^l \theta^j}(e^t,\theta) \right ] \partial_{\theta^l \theta^j}\tilde w_i^{\xi_i}(t,\theta), $$

and,

$$ D_i=\left [ \dfrac{1}{ \sqrt {| {\tilde g}^k|}(e^{t^{\xi_i}},\theta )  }           \partial_{\theta^l}[{\tilde g }^{\theta^l \theta^j}(e^{t^{\xi_i}},\theta)\sqrt {| {\tilde g}^k|}(e^{t^{\xi_i}},\theta)  ] -\dfrac{1}{ \sqrt {| {\tilde g}^k|}(e^t,\theta) } \partial_{\theta^l} [{\tilde g }^{\theta^l \theta^j}(e^t,\theta)\sqrt {| {\tilde g}^k|}(e^t,\theta) ] \right ] \partial_{\theta^j} \tilde w_i^{\xi_i}(t,\theta), $$

we deduce,

$$ A_i \leq C_k|e^{2t}-e^{2t^{\xi_i}}|\left [ |\nabla_{\theta} \tilde w_i^{\xi_i}| + |\nabla_{\theta}^2(\tilde w_i^{\xi_i})| \right ], $$

We take $ C=\max \{ C_i, 1 \leq i\leq q \} $ and if we use $ (***1) $, we obtain proposition 4.

We have:

$$ \dfrac{\partial_{\theta_j}w_i^{\lambda }(t,\theta)}{w_i^{\lambda }}=\dfrac{e^{(n-2)[(\lambda -\lambda_i)+(\xi_i-t)]/2} e^{[(\lambda -\lambda_i)+(\xi_i-t)]}(\partial_{\theta_j} v_i)(e^{[(\lambda -\lambda_i)+(\lambda -t)]}\theta)}{e^{(n-2)[(\lambda-\lambda_i)+(\lambda-t)]/2}v_i[e^{(\lambda-\lambda_i)+(\lambda - t)}\theta]} \leq {\bar C_i}, $$

 Also:
$$ \dfrac{\partial_{\theta_j,\theta_l}w_i^{\lambda }(t,\theta)}{w_i^{\lambda }}=\dfrac{e^{(n-2)[(\lambda -\lambda_i)+(\xi_i-t)]/2} e^{2[(\lambda -\lambda_i)+(\xi_i-t)]}(\partial_{\theta_j,\theta_l} v_i)(e^{[(\lambda -\lambda_i)+(\lambda -t)]}\theta)}{e^{(n-2)[(\lambda-\lambda_i)+(\lambda-t)]/2}v_i[e^{(\lambda-\lambda_i)+(\lambda - t)}\theta]} \leq {\bar C_i}. $$

where $ \bar C_i $ tends to $ 0 $ and does not depend on $ \lambda $.

We have,

$$ c(y_i,t,\theta)=\left ( \dfrac{n-2}{2} \right )^2+ \dfrac{n-2}{2} \partial_t a + R_g e^{2t}, \qquad (\alpha_1) $$ 

$$ b_2(t,\theta)=\partial_{tt} (\sqrt {b_1})=\dfrac{1}{2 \sqrt { b_1}}\partial_{tt}b_1-\dfrac{1}{4(b_1)^{3/2}}(\partial_t b_1)^2 ,\qquad (\alpha_2) $$ 

$$ c_2=[\dfrac{1}{\sqrt {b_1}} \Delta_{g_{y_i, e^t, {}_{{\mathbb S}_{n-1}}}}(\sqrt{b_1}) + |\nabla_{\theta} \log (\sqrt {b_1})|^2], \qquad (\alpha_3) $$

Then,

$$ \partial_{t}c(y_i,t,\theta)=\dfrac{(n-2)}{2}\partial_{tt}a, $$

by proposition 1,

$$ |\partial_tc_2|+|\partial_t b_1|+|\partial_t b_2|+|\partial_t c|\leq K_1e^{2t}. $$

We have:

$$ w_i(2\xi_i-t,\theta)=w_i[(\xi_i-t+\xi_i-\lambda_i-2)+(\lambda_i+2)] , $$

Thus,
 $$ w_i(2\xi_i-t,\theta)=e^{[(n-2)(\xi_i-t+\xi_i-\lambda_i-2)]/2}e^{n-2}v_i[\theta e^2e^{(\xi_i-t)+(\xi_i-\lambda_i-2)}]
\leq 2^{(n-2)/2}e^{n-2}=\bar c. $$

We set $ \delta = \dfrac{(n+2)-(n-2)\alpha}{2} $.

The left right side are denoted $ Z_1 $ et $ Z_2 $, we can write:

$$ Z_1=(\bar V_i^{\xi_i }-\bar V_i)(\tilde w_i^{\xi_i })^{N-1}+\bar
V_i[(\tilde w_i^{\xi_i })^{N-1}-{\tilde w_i}^{N-1}],$$

and,

$$ Z_2= e^{\delta
  t}[(\tilde w_i^{\xi_i })^{\alpha }-(\tilde w_i)^{\alpha
  }]+(\tilde w_i^{\xi_i})^{\alpha }(e^{\delta t^{\xi_i }}-e^{\delta t} ). $$

We can write the part with nonlinear terms as:

$  (\tilde w_i^{\xi_i})^{\alpha}[ (A\,  {w_i^{\xi_i}}^{N-1-\alpha}+ B)  \,  (e^t-e^{ t^{\xi_i }})+ c\,\,  (e^{\delta t^{\xi_i
    }}-e^{\delta t} ) ].
$
\smallskip

Because $ w_i^{\xi_i} \leq \bar c $, we have:

\smallskip

$$ -\bar Z_i((\tilde w_i)^{\xi_i}-(\tilde w_i))\leq  (\tilde w_i^{\xi_i})^{\alpha } [(A {\bar c}
^{N-1-\alpha}+ B)  \,  (e^t-e^{ t^{\xi_i }})+ c\,\,  (e^{\delta t^{\xi_i
    }}-e^{\delta t} ) ] + (\tilde w_i)^{\xi_i}(\mu/2)(e^{2t^{\xi_i
    }}-e^{2t} ) ] $$

\smallskip

But $ \alpha \in ] \dfrac{n}{n-2}, \dfrac{n+2}{n-2}[ $, $
\delta=\dfrac{n+2-(n-2)\alpha}{2} \in ]0,1[ $. 

\smallskip

We obtain for $ t \leq  t_0 <0 $:

$$ e^t \leq   e^{(1-\delta
  )t_0} e ^{\delta t} ,\,\,\, {\rm pour \,\, tout } \,\,\, t\leq t_0 .$$

and, $ t^{\xi_i}\leq t $ $(\xi_i \leq t )$, we integrate:

$$ (e^{\delta t^{\xi_i
    }}-e^{\delta t} ) \leq \dfrac{ \delta} { e^{(1-\delta
  )t_0} }\,  ( e^{ t^{\xi_i }}-e^t). $$

Finaly:

$$ -\bar Z_i((\tilde w_i)^{\xi_i}-(\tilde w_i)) \leq (\tilde w_i^{\xi_i})^{\alpha}[-\dfrac{\delta \, c }{ e^{(1-\delta
  )t_0} }+ A \,{\bar c}^{N-1-\alpha}+B]( e^t-e^{ t^{\xi_i }})+(\tilde w_i)^{\xi_i}(\mu/2)(e^{2t^{\xi_i
    }}-e^{2t}) $$

We apply  proposition 3. We take $ t_i=\log \sqrt {l_i} $ with $ l_i $ like in proposition 2. The fact $ \sqrt {l_i} [u_i(y_i)]^{2/(n-2)} \to +\infty $ ( see proposition 2), implies $ t_i=\log \sqrt {l_i} > \dfrac{2}{n-2} \log u_i(y_i) + 2 =\lambda_i+2 $. Finaly, we can work  on $ ]-\infty, t_i] $.

\bigskip

We define $ \xi_i $ by:

$$ \xi_i=\sup \{ \lambda \leq \lambda_i+2, \,\, \tilde w_i(2\lambda-t,\theta)-\tilde w_i(t,\theta)\leq 0\,\, {\rm on} \,\, [\lambda,t_i]\times {\mathbb S}_{n-1} \}. $$ 

If we use proposition 4 and the similar technics that in [B2] we can deduce by Hopf maximum principle,

$$ \min_{{\mathbb S}_{n-1}} \tilde w_i(t_i,\theta) \leq \max_{{\mathbb S}_{n-1}} \tilde w_i(2\xi_i-t_i,\theta), $$

which implies,

$$ {l_i}^{(n-2)/2} u_i(y_i) \times \min_M u_i \leq c. $$

It is in contradiction with  proposition 2.

\bigskip

Then we have,

$$ \sup_K u \times \inf_M u \leq c=c(\alpha, a, b, A, K, M, g, n). $$

\end{document}